\documentclass[12pt]{article}
\usepackage{amssymb}
\usepackage{amsmath}
\usepackage{amsthm}
\usepackage{mathscinet}
\usepackage[showonlyrefs]{mathtools}
\mathtoolsset{showonlyrefs}
\def\R{\mathbb{R}}
\def\N{\mathbb{N}}
\def\C{\mathbb{C}}
\def\arg{\operatorname{arg}}
\def\re{\operatorname{Re}}
\def\sgn{\operatorname{sgn}}
\newtheorem*{theorem}{Theorem}
\newtheorem*{lemma}{Lemma}
\theoremstyle{remark}

\newtheorem*{ack}{Acknowledgment}
\newtheorem{example}{Example}
\begin{document}
\title{On the number of solutions of some transcendental equations}
\author{Walter Bergweiler and Alexandre Eremenko\thanks{Supported by NSF
grant DMS-1665115.}}
\date{}
\maketitle
\emph{\centering{Dedicated to Dima Khavinson on the occasion of his 60th birthday}}
\begin{abstract}
We give upper and lower bounds for the number of solutions of the equation
$p(z)\log|z|+q(z)=0$ with polynomials $p$ and~$q$.
\end{abstract}

\section{Introduction and main result}\label{intro}
Holomorphic functions are sense-preserving. 
This allows, for a holomorphic function $f$ and $c\in\C$, to estimate
the number of solutions of the equation $f(z)=c$ from
above by the topological degree. This method does not work when $f$ is just smooth,
or real analytic, unless $f$ is sense-preserving. For the equation
\begin{equation}\label{1c}
\overline{z}=f(z),
\end{equation}
where $f$ is holomorphic,
a remarkable argument combining topological degree considerations 
with Fatou's theorem from holomorphic dynamics was invented by
Khavinson and \'Swi\polhk atek~\cite{KS}.
In this paper $f$ was a polynomial; later the argument was extended to 
rational $f$ by Khavinson and Neumann~\cite{KN}.
The latter result found an important
and unexpected application in astronomy. 
For transcendental meromorphic $f$ the equation~\eqref{1c} was
considered in \cite{BE2,BE,KL}, motivated by certain applications.
For a description of the method initiated in~\cite{KS} and its applications to
astronomy we also refer to the survey~\cite{K}.

This paper is a part of our efforts to understand the scope of
applicability of the method.
The following question was asked on Math Overflow~\cite{MO}.
Let $p$ and $q$ be coprime polynomials of degrees $m$ and~$n$, respectively,
with at least one of the polynomials non-constant.
How many solutions can the equation
\begin{equation}\label{1}
p(z)\log|z|+q(z)=0
\end{equation}
have?

\begin{theorem}
The number $N$ of solutions of equation~\eqref{1} satisfies
\begin{equation}\label{1b}
\max\{ m,n\}\leq N\leq 3\max\{ m,n\}+2m.
\end{equation}
\end{theorem}

The proof of the upper bound, given in section \ref{proof}, 
combines the computation of a topological degree
with Fatou's theorem as in the paper~\cite{KS}
mentioned above.
The difference of our argument in comparison
with previous applications of the method is that we transform~\eqref{1}
to an equation with infinitely many solutions, but it is
still possible to obtain the desired estimate.

The computation of the topological degree also yields the lower bound,
but only if solutions are counted with multiplicities.
In order to obtain a lower bound for the number of \emph{distinct} solutions
we study the curves where the rational function $q/p$ is real.

In section \ref{examples} we give examples to show that the estimate is sharp, 
at least for many values of $m$ and~$n$.

\begin{ack}
We thank the referee for helpful suggestions.
\end{ack}

\section{Proof of the theorem}\label{proof}
We put
\begin{equation}\label{f}
f(z)=\frac{2q(z)}{p(z)}
\end{equation}
and rewrite our equation~\eqref{1} as 
\begin{equation}\label{1p}
g(z):=\log|z|^2+f(z)=\log|z|^2+\frac{2q(z)}{p(z)}=0.
\end{equation}
The function $g$ is a continuous map of the Riemann sphere $\overline{\C}$ into itself,
satisfying $g(0)=g(\infty)=\infty$.

We recall the definition of the topological (or Brouwer) degree; see~\cite[Chapter~II, \S 2]{R}
or~\cite[\S 5]{Milnor1965}.
A value $w$ is \emph{regular} for $g$ if 
for all solutions $z$ of the equation $g(z)=w$ the map $g$
is continuously differentiable near $z$ and
the Jacobian determinant $J_g(z)$ does not vanish.
Then
\[
\deg g =\sum_{z\in g^{-1}(w)} \sgn J_g(z)
\]
is the \emph{topological degree} of~$g$.
This definition does not depend on~$w$.
(We note that e.g.\ in~\cite[\S\S 1--2]{Deimling} the topological degree is introduced for
functions on bounded domains, but this could be achieved by 
considering $g$ as map from $\{z\in\C\colon |f(z)|<R\}$ onto $\{z\in\C\colon |z|<R\}$
for some large~$R$.)

Taking $w=ir$ with large real $r$ we find $\max\{n-m,0\}$ 
preimages near $\infty$ and $m$ preimages near the poles. Since
\begin{equation}\label{3}
\begin{aligned}
J_g(z)&=|g_z(z)|^2-|g_{\overline{z}}(z)|^2
\\ &
=\left|\frac{1}{z}+f'(z)\right|^2-\frac{1}{|\overline{z}|^2}
\\ &
=\frac{1}{|z|^2}\left(|1+zf'(z)|^2-1\right)
\end{aligned}
\end{equation}
and since $|f'(z)|$ tends to $\infty$ as $z$ tends to a pole
and, if $n>m$, is bounded away from $0$ as $z$ tends to $\infty$,
we see that $J_g(z)>0$
at all these preimages, provided $r$ has been chosen sufficiently large.
So with 
\[
d:=\max\{ m,n\}=m+\max\{n-m,0\}
\]
we have $\deg g=d$.

For $w\in\C$ we denote by $N_w$
the number of solutions of 
\begin{equation}\label{2a}
g(z)=w
\end{equation}
so that $N=N_0$.
Suppose first that $w$ is a regular value of~$g$. 
We denote by $N^+$ and $N^-$ the numbers of solutions of~\eqref{2a} 
where $J_g(z)$ is positive  and negative, respectively.
Then
\begin{equation}\label{2}
N^+-N^-=d
\end{equation}
by the definition of the topological degree.

We put
\begin{equation}\label{4a}
h(z)=\frac{e^{-f(z)+w}}{z}
\end{equation}
and note that if $z$ satisfies~\eqref{2a}, then $z$ also satisfies
\begin{equation}\label{4}
k(z):=h(z)-\overline{z}=0.
\end{equation}
Note that the set of solutions of~\eqref{2a} is, in general, not equal to but only
{\em contained} in the set of solutions of~\eqref{4}.
The equation~\eqref{4} can have infinitely many solutions;
for example this is the case for the equation $\overline{z}=e^z/z$.

Since
\begin{equation}\label{h'}
h'(z)=-\frac{e^{-f(z)+w}}{z^2}(1+zf'(z))
\end{equation}
the Jacobian of $k$ is given by
\begin{equation}\label{second0}
J_k(z)=|h'(z)|^2-1=\frac{|e^{-2f(z)+2w}|}{|z|^4}|1+zf'(z)|^2-1.
\end{equation}
If $z$ is a solution of~\eqref{2a}, then $|z|^2=\exp(-f(z)+w)$ and thus
\begin{equation}\label{second}
J_k(z)=|1+zf'(z)|^2-1.
\end{equation}
We deduce from~\eqref{3} and~\eqref{second} that
the Jacobians $J_g(z)$ and $J_k(z)$ have the same sign for $z$ satisfying~\eqref{2a}.

Thus 
\begin{equation}\label{5}
N^-\leq n^-,
\end{equation}
where $n^-$ is the number of solutions of $k(z)=h(z)-\overline{z}=0$ with negative Jacobian. For these solutions we have
$|h'(z)|<1$, so they are exactly the {\em attracting fixed points} of
the antiholomorphic function $\overline{h(z)}$.

As already mentioned in the introduction, we will use Fatou's theorem from
complex dynamics. This theorem relates attracting fixed points to singular values.
To define singular values, we note that 
the essential singularities of $h$ are the poles of $f$.
We consider $h$ as a map from $\overline{\C}\backslash f^{-1}(\infty)$ to $\overline{\C}$.
If $h$ is not locally injective at a point $c\in \overline{\C}\backslash f^{-1}(\infty)$,
then $c$ is called a \emph{critical point} and $h(c)$ is called a \emph{critical value}.
In general, the set of critical points consists of the zeros of the derivative and
the multiple poles, but since our map $h$ has no multiple poles, we only have to consider
the zeros of~$h'$.
A value $a \in \overline{\C}$ is called an \emph{asymptotic value}
of $h$ if there is a path
$\gamma\colon [0,1)\to\overline{\C}\backslash f^{-1}(\infty)$ such that
$\gamma(t)$ tends to one of the essential singularies of $h$
as $t\to 1$ while $h(\gamma(t))\to a$ as $t\to 1$.
The singularities of the inverse function of $h$, or 
\emph{singular values} for short, are the critical and asymptotic values of~$h$.

They play an important role in complex dynamics.
The generalized Fatou theorem says that the basin of attration of an attracting
fixed point of a holomorphic (or antiholomorphic) function contains a singular value.
In particular,
the number of attracting fixed points of a holomorphic (or
antiholomorphic) function does not exceed the number of singular
values; see~\cite[Lemma~8.5]{Milnor2006} for rational functions, 
\cite[Lemma 10~(i)]{BDH} for functions which are meromorphic in $\overline{\C}$ except
for a compact, totally disconnected set (and thus in particular our
function~$h$), and \cite[p.~2914]{BE} for a version for self-maps of a Riemann surface
(which also applies to our function~$h$).

The number of singular values of $h$ is easy to estimate.
By~\eqref{h'},  the critical points of $h$  are the zeros of $1+zf'(z)$ in $\C$,
so there are at most $\max\{ m+n,2m\}=d+m$ of them. The asymptotic values
of $h$ can be only $0$ and $\infty$.
If $f(0)\neq\infty$, then $h(0)=\infty$. Otherwise, $0$ is an essential singularity of~$h$.
Moreover, if $f(\infty)\neq\infty$, then $h(\infty)=0$, while $\infty$
is an essential singularity of~$h$ if $f(\infty)=\infty$.
In any case we see that $0$ and $\infty$ either form a periodic cycle of period~$2$ 
for~$h$,
or they are essential singularities or mapped to essential singularities of~$h$.
In any case, they do not contribute to the count of attracting fixed points.
Thus 
$$N^-\leq n^- \leq d+m.$$
Combining this with~\eqref{2} we find that the number $N_w$ of solutions of~\eqref{2a} 
satisfies
\begin{equation}\label{2c}
N_w =N^+ +N^-= 2N^- +d 
\leq 2(d+m) +d=3d+2m.
\end{equation}
This proves the upper estimate in~\eqref{1b} if $0$ is a regular value of~$g$.

To deal with the case that $0$ is not regular we use the following lemma proved
in~\cite[Proposition~3]{BE2}.

\begin{lemma}
Let $D$ be a region in $\C$ and let $g\colon D\to\C$ be harmonic.
Suppose that there exists $M\in\N$ such that every $w\in\C$
has at most $M$ preimages under~$g$.
Then the set of points which have $M$ preimages is open.
\end{lemma}

We show that our function $g$ satisfies the hypothesis of this lemma
for a suitable domain~$D$. In order to do this we note that if $z$ satisfies~\eqref{2a},
then $z$ is a fixed point of the function 
\begin{equation}\label{2b}
\zeta\mapsto \overline{h\!\left(\overline{h(\zeta)}\right)}.
\end{equation}
This function is holomorphic in $\C$ except for singularities at $0$ and the poles of~$f$.
So the solutions of~\eqref{2a} form a discrete set.
Since the solutions of~\eqref{2a} do not accumulate at~$0$, $\infty$ or a pole of~$f$, we 
conclude that~\eqref{2a} has only finitely many solutions, for each~$w\in\C$.
We thus have $N_w<\infty$ also if $w$ is not regular; that is, for each $w\in\C$ the
function $g$ has only finitely many $w$-points. 
In order to apply the lemma we still have to show that
the number of $w$-points is uniformly bounded by some $M\in\N$,
at least after restricting to a suitable domain~$D$.

We denote by $N_w(D)$ the number of $w$-points of $g$ in a domain~$D$.
We choose a bounded domain $D$ containing all solutions of the equation $g(z)=0$ such that 
the closure of $D$ does not contain $0$ or a pole of~$f$. By the choice of $D$ we then
have $N=N_0(D)$.
If $\zeta\in D$ is such that $J_g(\zeta)\neq 0$,
then $\zeta$ clearly has a neighborhood $U_\zeta$ such that $N_w(U_\zeta)\leq 1$ for all $w\in\C$.
Moreover, it follows from results of Lyzzaik~\cite[Theorems~5.1 and~6.1]{Lyzzaik1992} that
if $\zeta\in D$ with $J_g(\zeta)= 0$, then there exist a neighborhood $U_\zeta$ of
$\zeta$ and $M_\zeta\in \N$ such that $N_w(U_\zeta)\leq M_\zeta$ for all $w\in\C$.
(The results of Lyzzaik give precise information about the value~$M_\zeta$, but this is
irrelevant for our purposes.)
Since $D$ can be covered by finitely many neighborhoods $U_\zeta$, we deduce that
there exists $M\in\N$ such that $N_w(D)\leq M$ for all $w\in\C$.
We may assume that $M$ has been chosen minimal. Then the set of all $w\in \C$ with
$N_w(D)=M$ is a non-empty open subset of~$\C$ by the lemma. This implies 
that there exists a regular value $w$ with $N_w(D)=M$.
Combining this with~\eqref{2c} we thus have
\begin{equation}\label{2d}
N=N_0(D)\leq M=N_w(D)\leq N_w
\leq 3d+2m.
\end{equation}
This shows that the upper estimate in~\eqref{1b} also holds if $0$ is not a regular value.

To prove the lower estimate in~\eqref{1b} we put 
\begin{equation}\label{5a}
F(z)=\frac{q(z)}{p(z)} =\frac12 f(z).
\end{equation}
So $F$ is a rational function of degree~$d$.
If $F$ has no real critical values, the preimage of $\R$ under $F$ is a union of 
$d$ disjoint curves in~$\overline{\C}$. The start and end point of such a curve are (not 
necessarily distinct) poles. If $F(\infty)$ is finite and real, then at least one and possibly
several of these curves pass through~$\infty$.

If $F$ has real critical values, we consider these curves for the function
$F-i\varepsilon$ instead of~$F$, for some small positive~$\varepsilon$.
Taking the limit as $\varepsilon\to 0$ we find that $F^{-1}(\R)$ 
is still the union of $d$ curves $\gamma_1,\dots,\gamma_d$, with
each $\gamma_j$ starting and ending at a pole of~$F$,
but now these curves are not disjoint anymore.

Indeed, let $w_0\in\R$ be a critical value of~$F$, say $w_0=F(z_0)$ where $z_0\in\C$ with $F'(z_0)=0$.
Let $L$ be the multiplicity of the $w_0$-point $z_0$; that is,
$L=\min\{k\in\N\colon F^{(k)}(z_0)\neq 0\}$.
Then there exist $L$ curves passing through $z_0$, and we may assume that the curves are numbered
so that this is the case for the curves $\gamma_1,\dots,\gamma_L$.
Choosing parametrizations $\gamma_j\colon I_j\to\overline{\C}$ with 
intervals  $I_j$ we thus have $\gamma_{j}(t_{j})=z_0$ for some $t_{j}\in I_j$.
We may assume that the parametrizations $\gamma_j$ are chosen 
such that $F(\gamma_j(t))$ increases with~$t$.
The directions of the curve $\gamma_{j}$ at the point $z_0$ are
given by the one-sided derivatives $\gamma_{j}'(t_{j}^\pm)$ of $\gamma_{j}$ at~$t_{j}$.
The left and right derivative are related by
\begin{equation}\label{5b}
\arg\gamma_{j}'(t_{j}^+)=\arg\gamma_{j}'(t_{j}^-)+\frac{\pi}{L}-\pi.
\end{equation}
Moreover, for a suitable permutation $\sigma\in S_L$ we have
\begin{equation}\label{5c}
\arg\gamma_{j}'(t_{j}^+)=\frac{2 \pi \sigma(j)}{L} -\alpha
\end{equation}
where $\alpha=\arg F^{(L)}(z_0)$.

We now consider the function 
\begin{equation}\label{5c1}
G_j\colon I_j\to \R, \quad G_j(t)=F(\gamma_j(t))+\log|\gamma_j(t)|.
\end{equation}
Noting that there are poles $p_j^+$ and $p_j^-$ such that $\gamma_j(t)\to p_j^+$ as $t\to \sup I_j$
while $\gamma_j(t)\to p_j^-$ as $t\to \inf I_j$
we can deduce that $G_j(t)\to \pm\infty$ as $t\to \sup I_j$ or $t\to\inf I_j$, respectively.
This is clear if $p_j^-\neq 0$ and $p_j^+\neq \infty$, but it also follows if
$p_j^-= 0$ or $p_j^+= \infty$, since then $F(\gamma_j(t))$ tends to $\pm\infty$
faster than $\log|\gamma_j(t)|$.

Thus there exists $s_j\in I_j$ such that $G_j$ changes its sign from $-$ to $+$
at~$s_j$; that is, there exists $\delta>0$ such that $G_j(s)<0$ for $s_j-\delta<s<s_j$ 
while $G_j(s)>0$ for $s_j<s<s_j+\delta$. 
It follows that
\begin{equation}\label{5d}
G_j'(s_j^-)\geq 0 \quad\text{and}\quad G_j'(s_j^+)\geq 0 .
\end{equation}
If $\gamma_j$ passes through~$\infty$, which can happen only if $F(\infty)$ is finite and real, then 
$F(z)-\log|z|$ is negative for all $z$ on this curve of sufficiently large modulus.
This implies that $\gamma_j(s_j)\in\C$ and hence
$z=\gamma_j(s_j)$ is a solution of our equation~\eqref{1}.
If all the points $\gamma_j(s_j)$ are distinct we thus have $d$ solutions.
This is clearly the case if none of the points $\gamma_j(s_j)$ is a critical point.

Suppose now that $z_0=\gamma_j(s_j)$ is a critical point for some~$j$.
Using the notation above we thus have $j\in\{1,\dots,L\}$ and $s_j=t_{j}$.

Noting that
\begin{equation}\label{5e}
\frac{d}{dt} \log|\gamma_{j}(t)|=\re\frac{\gamma_j'(t)}{\gamma_j(t)} 
\end{equation}
and $F'(z_0)=0$ we then have
\begin{equation}\label{5f}
\begin{aligned}
G_j'(t_{j}^\pm)
& = F'(\gamma_{j}(t_{j})) \gamma_{j}'(t_{j}^\pm)+
\re\frac{\gamma_{j}'(t_{j}^\pm)}{\gamma_{j}(t_{j})} 
\\ &
= F'(z_0) \gamma_{j}'(t_{j}^\pm)+
\re\frac{\gamma_{j}'(t_{j}^\pm)}{z_0} 
\\ &
=\re\frac{\gamma_{j}'(t_{j}^\pm)}{z_0}.
\end{aligned}
\end{equation}
Put $\beta=\arg z_0$. In view of~\eqref{5d} the last equation yields that
\begin{equation}\label{5h}
\cos\!\left(\arg \gamma_{j}'(t_{j}^+)-\beta\right)\geq 0
 \quad\text{and}\quad 
\cos\!\left(\arg \gamma_{j}'(t_{j}^-)-\beta\right)\geq 0 .
\end{equation}
Since
\begin{equation}\label{5i}
\begin{aligned}
\cos\!\left(\arg \gamma_{j}'(t_{j}^-)-\beta\right)
&=
\cos\!\left(\arg \gamma_{j}'(t_{j}^+)-\frac{\pi}{L}+\pi-\beta\right)
\\&
=
-\cos\!\left(\arg \gamma_{j}'(t_{j}^+)-\frac{\pi}{L}-\beta\right)
\end{aligned}
\end{equation}
by~\eqref{5b} we deduce from~\eqref{5c} and~\eqref{5h} with $\theta=\alpha+\beta$ that
\begin{equation}\label{5j}
\cos\!\left(\frac{2\pi\sigma(j)}{L} -\theta\right)\geq 0
 \quad\text{and}\quad 
\cos\!\left(\frac{2\pi\sigma(j)}{L} -\frac{\pi}{L}-\theta\right)\leq 0 .
\end{equation}
Since in an interval of length $2\pi$ there is only point where the 
cosine changes its sign from $-$ to $+$ there exists at most
one value $\sigma(j)\in\{1,\dots,L\}$ that satisfies~\eqref{5j}.
We conclude that if $z_0$ is a critical point of~$F$, then
$z_0=\gamma_j(s_j)$ for at most one value of~$j$.
Altogether we see that the points $\gamma_j(s_j)$ are all distinct so 
that our equation has at least $d=\max\{m,n\}$ solutions.
This completes the proof of the theorem.

\section{Examples}\label{examples}
We give several examples to show that the estimates in our theorem are best 
possible.
More specifically, Examples~\ref{ex1} and~\ref{ex2} show that the upper bound
is sharp if $m=0$ or $n=0$. Example~\ref{ex3} deals with the case $n\leq m$, thus
generalizing Example~\ref{ex2}. Examples~\ref{ex4} and~\ref{ex5} show that the upper bound
is sharp if $n=2m$ or $n=3m$.
Finally, Example~\ref{ex6} shows that the lower bound is sharp for all $m$ and~$n$.

\begin{example} \label{ex1}
For  $p(z)=1$ and $q(z)=2\log 2 \cdot (1-z)$ the equation~\eqref{1} has the positive 
solutions $1/2$ and~$1$, and there is one negative solution $\xi$ by the 
intermediate value theorem. Computation shows that $\xi\approx -0.191666$.
This shows that the upper bound in the theorem is sharp for $m=0$ and $n=1$.
Considering
\[
q_n(z):=\frac1n q(z^n) =\frac{2\log 2}{n}(1-z^n)
\]
with $n\geq 2$ instead of $q$ we see that the upper bound is sharp
for $m=0$ and arbitrary $n\in\N$.

Indeed, for any $n$-th root of unity $\omega$
the equation $\log|z|+q_n(z)=0$ has the solutions $\omega$, $\omega/\sqrt[n]{2}$ 
and $\omega\sqrt[n]{|\xi|} e^{i\pi/n}$ so that there are $3n$ solutions 
altogether;
that is, the equation 
\begin{equation}\label{dwa2}
\log|z|+\frac{2\log 2}{n}(1-z^n)=0
\end{equation}
has $3n$ solutions.
\end{example}

\begin{example} \label{ex2}
For $p(z)=8z+1$ and $q(z)=6\log 2$
the equation~\eqref{1} has the three positive solutions $1/16$, $1/8$ and $1/4$,
and two negative solutions $\xi_{1,2}$ by the intermediate value theorem.
The numerical values are
$\xi_1\approx -1.471293$ and 
$\xi_2\approx -0.0106199$.
Similarly as in the previous example we see by 
considering $p_m(z)=m\; p(z^m)$ with $m\geq 2$ instead of $p$ 
that the upper bound  in our result is sharp for $n=0$ and arbitrary $m\in\N$;
that is, the equation 
\begin{equation}\label{dwa}
\log|z|+\frac{6\log 2}{m(8z^m+1)}=0
\end{equation}
has $5m$ solutions.
\end{example}

\begin{example} \label{ex3}
The previous example can be perturbed as follows. Choose a polynomial
$q$ of degree $n\leq m$ which is close to $1$ on a compact set
containing all $5m$ solutions of~\eqref{dwa}. As all solutions
of~\eqref{dwa} are non-degenerate, the inverse function theorem will
guarantee that the number of solutions of
\begin{equation}\label{dwa3}
\log|z|+\frac{6\log 2\; q(z)}{m(8z^m+1)}=0
\end{equation}
is at least $5m$ when $q$ is sufficiently close to~$1$.
This shows that the upper estimate in the theorem is best possible
for all $n\leq m$.
 
An explicit example with $m=n$ is
\begin{equation}\label{tri}
\log|z|=3\log 2\frac{z^n-1}{n(z^n+1)}.
\end{equation}
When $n=1$ this equation has $5$ real solutions: the positive solutions~$1$, $2$ and $1/2$,
as well as two negative solutions by the intermediate value theorem, which can be computed to
be 
$\xi_1\approx -11.770347$ and 
$\xi_2\approx -0.0849592$.
Making the change of the variable $z\mapsto z^n$ we see that ~\eqref{3}
has $5n$ solutions.
\end{example}

\begin{example} \label{ex4}
Take $a=0.015$, and consider the equation
\begin{equation}\label{chetyre}
\log|z|=3\log2\cdot (1-a(z-1))\frac{z-1}{z+1}.
\end{equation}
This is a small perturbation of~\eqref{tri} with $n=1$.
Again $z=1$ is clearly a solution
and one can check 
that it has $4$ further real solutions near the solutions of~\eqref{tri}.
Moreover, the intermediate value theorem yields that it has one more
negative solution. The numerical values of these  $6$ real zeros $\xi_1,\dots,\xi_6$
are at $\xi_1 \approx -58.249375$,
$\xi_2 \approx -20.915701$,
$\xi_3 \approx -0.0826000$,
$\xi_4 \approx 0.466285$,
$\xi_5 = 1$ and $\xi_6 \approx 1.780021$.

Let $f$ be the right hand side of~\eqref{chetyre}.
Then $f$ has two real critical points
$x_1\approx -12.718930$
and
$x_2\approx 10.718930$
with critical values
$y_1\approx 2.935272$
and
$y_2\approx 1.473143$.

This shows that there is a curve $\gamma$
in the upper half-plane with endpoints
$x_1$ and $x_2$ on which $f$ is real.
As 
$\log|x_1|\approx 2.543091478<y_1$
and
$\log|x_2|\approx 2.372011>y_2$
we conclude that the equation~\eqref{chetyre} must have a solution
in the upper half-plane and, by symmetry,
another one in the lower half-plane. 
Numerically these two solutions are
$\xi_{7,8}\approx -5.705306\pm 10.732819 i$.
Altogether the total number
of solutions of~\eqref{chetyre} is thus~$8$.

Making the change of variable $z\mapsto z^m$, we obtain
an equation with $n=2m$ having $8m=3\cdot 2m+2m$
solutions. This shows that the upper estimate in the theorem  is exact when
$n=2m$.
\end{example}

\begin{example} \label{ex5}
This example is again a small perturbation of the previous example. 
As there we take $a=0.015$,
put $b=0.00185$ and and consider the equation
\begin{equation}\label{chetyre2}
\log|z|=3\log2\cdot (1-a(z-1))\cdot (1-b(z-1))\frac{z-1}{z+1}.
\end{equation}
The equation has $7$ real solutions, $6$ of which correspond to the 
solutions of~\eqref{chetyre}.
The numerical values are 
$\xi_1 \approx -198.8150$,
$\xi_2 \approx -176.4617$,
$\xi_3 \approx -17.8054$,
$\xi_4 \approx 0.08289$,
$\xi_5 \approx 0.4704$,
$\xi_6 = 1$ and $\xi_7 \approx 1.8020$.
Denoting by $f$ the right hand side of~\eqref{chetyre2} we see that $f$
has two critical points near those found in the previous example,
and  there is a curve connecting these points in the upper half-plane on which $f$ is real.
On this curve we then have a solution of~\eqref{chetyre2}.
Together with its complex conjugate this yields the two solutions
$\xi_{8,9}\approx 8.6167\pm 10.2654i$.

Moreover, $f$ has one critical point at $x_0\approx -234.2572$, 
and we have $f(x_0)<\log|x_0|$. This yields that there exists a curve
in the upper half-plane connecting $x_0$ with $\infty$ on which $f$ is real.
This curve then contains a solution of~\eqref{chetyre2}.
Together with its complex conjugate we obtain the solutions
$\xi_{10,11}\approx -234.2803\pm 43.6244i$.

Altogether we thus have $11$ solutions. 
The change of variable $z\mapsto z^m$ then yields 
an equation with $n=3m$ having $11m=3\cdot 3m+2m$
solutions. Thus the upper estimate in the theorem  is exact when $n=3m$.
\end{example}

\begin{example} \label{ex6}
Let $p$ and $q$ be polynomials of degrees $m$ and~$n$, respectively.
Suppose that $F(z):=q(z)/p(z)\in\C\backslash\R$ for $|z|\leq 1$.
If $F(\infty)\in\C\backslash\{0\}$, assume in addition that $F(\infty)\notin\R$.
It is clear that polynomials $p$ and $q$ with these properties exist.
In fact, if $p_0$ and $q_0$ are polynomials satisfying $p_0(0)\neq 0$ and $q_0(0)\neq 0$,
then there exists $\varphi\in\R$ such that $p(z)=e^{i\varphi}p_0(\delta z)$ and
$q(z)=e^{i\varphi}q_0(\delta z)$ have the above properties for all small positive~$\delta$.

We show that if $c$ is a large positive number, then the equation
\begin{equation}\label{1a}
c p(z)\log|z|+c q(z)=0
\end{equation}
has $\max\{m,n\}$  solutions.
This shows that the lower bound in our theorem is best possible.

As before we put $d=\max\{m,n\}$. For $1\leq j\leq d$ we choose the curves
$\gamma_j\colon I_j\to \overline{\C}$
as in the proof of the theorem.
Since $F(z)\in\C\backslash\R$ for $|z|\leq 1$ we find that the curves $\gamma_j$
are contained in $\{z\colon|z|>1\}\cup\{\infty\}$.
Since $F$ has no real critical values, we have
\begin{equation}\label{6a}
|F'(z)|>0\quad\text{if}\ z\in\C \ \text{and}\ F(z)\in\R.
\end{equation}
Suppose first that $j$ is such that the curve $\gamma_j$ does not pass through~$\infty$.
For a sufficiently large positive constant $c$ we then can have
\begin{equation}\label{6a1}
c|F'(\gamma_j(t))|> \frac{1}{|\gamma_j(t)|}
\quad\text{for}\ t\in I_j.
\end{equation}
Note that this also works if $\infty$ is an endpoint of $\gamma_j$ since then $F(\infty)=\infty$.
Let $G_{c,j}$ be defined as in~\eqref{5c1}, with $F$ replaced by $cF$; that is,
\begin{equation}\label{5c2}
G_j\colon I_j\to \R, \quad G_j(t)=c F(\gamma_j(t))+\log|\gamma_j(t)|.
\end{equation}
We deduce from~\eqref{6a1} that
\begin{equation}\label{6c}
G_{c,j}'(t)
 = cF'(\gamma_{j}(t)) \gamma_{j}'(t)+
\re\frac{\gamma_{j}'(t)}{\gamma_{j}(t)}
 \geq |\gamma_j'(t)|\left (c| F'(\gamma_{j}(t))|
-\frac{1}{|\gamma_{j}(t)|} \right)>0.
\end{equation}
Hence $G_{c,j}$ is increasing and thus $G_{c,j}$ has exactly one zero.

Suppose now that $\gamma_j$ passes through $\infty$, say $\gamma_j(t_j)=\infty$.
Noting that $F(\infty)=F(\gamma_j(t_j))\notin\R$ if $F(\infty)\in\C\backslash\{0\}$,
by our choice of $p$ and $q$, we see that $F(\infty)=0$. 
Next we note that as it approaches $\infty$, the curve $\gamma_j$ is asymptotic to a ray from the origin.
In fact, it is not difficult to show that 
\begin{equation}\label{6e}
\re\frac{\gamma_{j}'(t)}{\gamma_{j}(t)}\sim \left|\frac{\gamma_{j}'(t)}{\gamma_{j}(t)}\right|
\quad\text{as}\ t\to t_j, \ t<t_j.
\end{equation}
In particular, there exists $s_j\in I_j$ with $s_j<t_j$ such that 
\begin{equation}\label{6e1}
\re\frac{\gamma_{j}'(t)}{\gamma_{j}(t)}>0
\quad\text{for}\  s_j\leq t<t_j.
\end{equation}
Since $F(\gamma_j(t))$ increases with $t$ we have
\begin{equation}\label{6e2}
\frac{d}{dt}F(\gamma_j(t))=F'(\gamma_{j}(t))\gamma_{j}'(t)\geq 0.
\end{equation}
The last two inequalities imply that 
\begin{equation}\label{6f}
G_{c,j}'(t)
>0
\quad\text{for}\  s_j\leq t<t_j.
\end{equation}

On the other hand, for $t\in I_j$ with $t\leq s_j$ we have
\begin{equation}\label{6g}
F(\gamma_j(t))\leq F(\gamma_j(t_j))<0.
\end{equation}
This implies that if $c$ is sufficiently large, then
\begin{equation}\label{6f1}
c F(\gamma_j(t))<-\log|\gamma_j(t)|
\quad\text{for}\  t\leq s_j.
\end{equation}
It follows from~\eqref{6f1} and~\eqref{6f} that 
$G_{c,j}(t)<0$ for $t\leq s_j$ 
and that $G_{c,j}$ is strictly increasing in $[s_j,t_j]$.
Moreover, $G_{c,j}(t)>0$ for $t>t_j$ since 
$F(\gamma_j(t))$ increases with $t$ and $F(\gamma_j(t_j))=0$ and since
$\gamma_j$ does not intersect $\{z\colon|z|\leq 1\}$.
Thus $G_{c,j}$ has exactly one zero also in this case.

Altogether we see that $G_{c,j}$ has exactly one zero for each $j\in\{1,\dots,d\}$.
Thus the equation $cF(z)+\log|z|=0$ has exactly $d$ zeros, from which the conclusion
follows.
\end{example}

\bigskip

\noindent Mathematisches Seminar\\
Christian-Albrechts-Universit\"at zu Kiel\\
Ludewig-Meyn-Str.\ 4\\
24098 Kiel\\
Germany

\bigskip

\noindent Department of Mathematics\\
Purdue University\\
West Lafayette, IN 47907\\
USA
\end{document}